\documentclass[12pt]{l4dc2022}
\usepackage[shortlabels]{enumitem}

\usepackage{amsmath}
\usepackage{mathtools}
\title[Sketching for system identification]{Learning Linear Models Using Distributed Iterative Hessian Sketching}
\usepackage{times}
\usepackage{algorithm}
\usepackage{todonotes}
\usepackage{algpseudocode}
\usepackage{hyperref}

\newcommand{\xls}{X^{\mathrm{LS}}}

\newcommand{\iid}{\stackrel{\mathclap{\text{\scriptsize{ \tiny i.i.d.}}}}{\sim}}
\newcommand{\R}{\mathbb{R}}

\newcommand{\sg}{\mathtt{subG}}

\renewcommand{\vec}{\mathrm{vec}}

\newcommand{\ra}{\text{Range}}

\author{%
 \Name{Han Wang} \Email{hw2786@columbia.edu}\\
 \addr Columbia University, New York, NY
 \AND
 \Name{James Anderson} \Email{james.anderson@columbia.edu}\\
 \addr Columbia University, New York, NY%
}

\begin{document}

\maketitle

\begin{abstract}

This work considers the problem of learning the Markov parameters of a linear system from observed data. Recent  non-asymptotic system identification results have characterized the sample complexity of this problem in the single and multi-rollout setting. In both instances, the number of samples required in order to obtain acceptable estimates can produce optimization problems with an intractably large number of decision variables for a second-order algorithm. We show that a randomized and distributed Newton algorithm based on Hessian-sketching can produce  $\epsilon$-optimal solutions and converges geometrically. Moreover, the algorithm is trivially parallelizable. Our results hold for a variety of sketching matrices and we illustrate the theory with numerical examples.

\end{abstract}


\section{Introduction}
Obtaining a dynamic model of a system or process is fundamental to most of science and engineering. As the systems we study become increasingly complex,  data-driven modeling has become the de facto framework for obtaining accurate models~\citep{BruK19}.  Fortunately, as systems become more interconnected and sensors become smaller and cheaper, there is no shortage of data to work with.  Indeed, the volume of data available can overwhelm the (often limited) computational resources at our disposal, forcing us to consider data versus  resource trade-offs~\citep{ChaJ13}.

There has recently been considerable interest in applying machine learning techniques to the problem of controlling a dynamical system. Two paradigms have emerged; \emph{model-based control}, in which a model is first learnt from data and then a classical controller is synthesized from the model. In the \emph{model-free} setting  the control action is learnt directly from data without ever constructing an explicit model, see for example~\cite{FazGHM18}. Our work is motivated by two observations; i) the asymptotic sample complexity of a model-based solution outperforms that of a model-free Least-Squares Temporal Difference Learning~\citep{Boy99} approach ~\cite{TuR19}; ii) the recent body of work characterizing the sample complexity of learning linear system models from data, shows that for systems with a large number of inputs and outputs, the resulting optimization problems are intractable as they require a large number of rollouts or long trajectory horizon lengths in order to produce accurate estimates~\cite{oymak2019non,zheng2020non,tsiamis2019finite,DeaMMRT20}. The goal of this work is to \emph{construct and solve approximations of these optimization problems} that are consistent with the sample complexity results, provide provably good solutions, \emph{and} do so in an algorithmically tractable manner. Our approach is based on the concept of ``sketching''~\citep{Drim16}.  Broadly speaking, a ``sketch'' is an approximation of a large matrix by a smaller or more ``simple'' matrix. For a sketch to be useful, it must retain certain properties of the original matrix that allow it to be used for computation in place of the original. What is perhaps surprising, is that \emph{randomization} is the enabling force used to construct sketches~\citep{Woo14, MarT20,Mah11}. Moreover, numerical linear algebra routines based on randomization (and sketching) can outperform their deterministic counterparts~\citep{AvrMT10}, and are more suited to distributed computing architectures.

\subsection{Problem Setting}
\textbf{Notation:}\ Given a matrix $A \in \mathbb{R}^{m\times n},$ we use $\lVert A \rVert_F$ to denote its Frobenius norm. The multivariate normal distribution with
mean $\mu$ and covariance matrix $\Sigma$ is denoted by $\mathcal{N}\left(\mu, \Sigma \right).$ For two functions $f(x)$ and $g(x)$, the notation $f(x)=O(g(x))$ or $f(x) \lesssim g(x)$ implies that there exists a universal constant $C<\infty$ satisfying $f(x) \leq C g(x)$. For an event $\mathcal{X}$,  $\mathbb{P}(\mathcal{X})$ refers to its probability of occurrence.\\

Let us assume that we have a stable and minimal, linear time-invariant (LTI) system given by
\begin{equation}\label{eq:fully observed}
	\begin{aligned}
		x_{t+1} &=A x_{t}+B u_{t}+ w_{t}, \quad x_0=0 \\
		y_{t} &=C x_{t}+D u_{t}+v_{t}
	\end{aligned}
\end{equation}
where $x_{t} \in \mathbb{R}^{n}, ~y_{t} \in \mathbb{R}^{p},$ and  $u_{t} \in \mathbb{R}^{m}$ denote the system state, output, and input at time $t$, respectively, and $w_{t} \in \mathbb{R}^{n}, v_{t} \in \mathbb{R}^{p}$ denote the process and measurement noises. We assume that ${u}_{t} \sim \mathcal{N}\left(0, \sigma_{u}^{2} {I}_{m}\right), {w}_{t} \sim \mathcal{N}\left(0, \sigma_{w}^{2} {I}_{n}\right)$, and ${v}_{t} \sim \mathcal{N}\left(0, \sigma_{v}^{2} {I}_{p}\right).$  Our goal is to learn the system parameters $A,B,C$ and $D$ from a single input and output trajectory $\{y_t,u_t\}^{\bar{N}}_{t=1}$. We are particularly interested in the scenario of ``large'' $m$ and $p$, and where sample complexity results show that $\bar N$ must be huge in order to achieve accurate estimates. In this parameter regime many optimization methods are intractable and even routine matrix factorizations become problematic.

To achieve this goal, we  begin by estimating the first $T$ Markov parameters $G$, which are defined as: 
$$G=\left[\begin{array} {ccccc}D & CB & CAB & \cdots & C A^{T-2} B\end{array}\right] \in \mathbb{R}^{p \times mT}.$$ Then system matrices $A,B,C$ and $D$ can be realized via the Ho-Kalman algorithm (\cite{ho1966effective}). In recent work we applied similar ideas based on randomized methods to implement a stochastic version of the Ho-Kalman algorithm suitable for masive-scale problems~\citep{WanA21}. We thus narrow our attention in this work to the task of providing an estimate $\hat G$ of $G$.

As described in~\cite{oymak2019non}, we first generate a trajectory of length $\bar N$. The trajectory is then spliced and written into two matrices corresponding to the control and output signal. Define $\bar{N} = T+N-1$ with $N \ge 1$, let
\begin{equation*}
Y=\left[\begin{array}{cccc}y_T & y_{T+1} & \cdots & y_{\bar{N}}\end{array} \right]^T \in \mathbb{R}^{N \times p} \text{~~~ and~~~ }
U=\left[\begin{array}{cccc}\bar{u}_T & \bar{u}_{T+1} & \cdots & \bar{u}_{\bar{N}}\end{array}\right]^T \in \mathbb{R}^{N \times mT},
\end{equation*}
with  $\bar{u}_i$ denoting $
\bar{u}_{t}=\left[
u_{i}^{T}, u_{i-1}^{T},\cdots,u_{i-T+1}^{T}
\right]^{T} \in \mathbb{R}^{mT}
$.
Then the Markov parameters $G$ can be learned by solving the following unconstrained least-squares problem:
\begin{equation}\label{eq:leastsq}
	\xls=\underset{{X} \in R^{mT \times p}}{\operatorname{argmin}}\|{Y}-{U} X\|_{F}^{2}.
\end{equation}
The estimate  is obtained from $\hat G =(\xls)^T$. \cite{oymak2019non} provided the following non-asymptotic sample complexity bound: 
\begin{equation}\label{eq:bound}
	\|\hat{G}-G\|_{F} \le \frac{\left(\sigma_{v}+\sigma_{e}\right) \sqrt{p}+\sigma_{w}\|F\|_{2}}{\sigma_{u}} \sqrt{\frac{T q \log ^{2}(T q) \log ^{2}(N q)}{N}}	
\end{equation}
holds with high probability, as long as $
N \gtrsim T q \log ^{2}(T q) \log ^{2}(N q)$, where $q = m+p +n$ is the aggregated system dimension. The matrix $
F=\left[\begin{array}{ccccc}
	0 & C & CA & \cdots & C A^{T-2}\end{array}
\right] \in \mathbb{R}^{p \times T n}
$ is the  concatenated matrix and $\sigma_e^2$ is the variance of the linearly transformed state at time $i-T+1.$  Interested readers can refer to~\cite{oymak2019non} for more details.

\subsection{Motivation}
From~\eqref{eq:bound}, it is clear that increasing the sample size $N$ can make the estimated Markov parameters more reliable. To achieve  better identification performance, $N$ needs to be  large, i.e., $N\gg T q \log ^{2}(T q) \log ^{2}(N q)$. In other words, we need to solve the least square problem described by Eq~\ref{eq:leastsq} with $N \gg mT$ (The number of rows is significantly larger than the number of columns).

To solve  problem~\eqref{eq:leastsq}, \cite{oymak2019non,zheng2020non,tsiamis2019finite,DeaMMRT20} adopted the pseudo-inverse method. However, the complexity of computing the pseudo-inverse method requires $O( (mT)^2N)$ flops, which is costly  for large  systems. Moreover, for truly huge-scale systems, the memory cost for storing the sample trajectories ($U$ and $Y$) will likely exceed the storage capacity of a single machine. Therefore, there is a strong desire to put forward a tractable algorithm, which can take advantage of modern distributed computing architectures. In this paper, we consider the setting where there are $r$ worker machines operating independently in parallel and a single central node that computes the averaged solution. No communication between  workers is permitted as it is likely that communication time dominates local computation time in the distributed algorithms. We only allow the communications between worker and the central node. Overall, we aim to provide a communication $\&$ computation-efficient algorithm to solve the large-scale system identification problems defined by Eq~\eqref{eq:leastsq}.

\subsection{Related work}
\textbf{System identification:}\ Estimating a linear dynamical system from input/output observations has a long history, which can date back to the 1960s. Prior to the 2000s, most identification methods for linear systems either focus on the prediction error approach~\cite{ljung1999system}
or  subspace methods~\cite{van2012subspace,verhaegen2007filtering}.
In contrast, with the advances in high-dimensional statistics~\cite{vershynin2018high}, contemporary research shifts from asymptotic
analysis with infinite data assumptions to finite time analysis and finite data rates. Over the past
several years, there have been significant advances in studying the finite sample
properties, when the system state is fully observed~\cite{simchowitz2018learning,sarkar2019near,faradonbeh2018finite}. When the system is partially observed, we can find the finite sample analysis in~\cite{oymak2019non,sarkar2019finite,simchowitz2019learning,tsiamis2019finite,lee2020non,zheng2020non,lee2020improved,lale2020logarithmic,kozdoba2019line}. However, there are only a few papers~\cite{sznaier2020control,reyhanian2021online} that consider the computational complexity and memory issues of system identification, which become prohibitively large and incompatible with on-board resources when system dimension increases. There is thus a great need to provide scalable algorithms which can efficiently solve the system identification problem.\\
\\
\textbf{Distributed optimization:}\ In recent years, a lot of effort has been devoted to designing distributed first-order methods~\citep{mahajan2013parallel,shamir2014distributed,lee2017distributed,fercoq2016optimization,liu2014asynchronous,necoara2016parallel,richtarik2016parallel,liu2020double}, which only rely on gradient information of the objective function. However, first-order methods  suffer from: (i) a dependence on a suitably defined condition number; (ii) spending more time on communication than on computation. To overcome these drawbacks, second-order methods have received more attention recently, since they enjoy superior convergence rates which are independent of the condition number and thereby require fewer rounds of communication to achieve high accuracy solutions.

 The trade-off is that  most second-order algorithms based on Newton's method require forming and then computing the inverse of  Hessian matrix at each iteration. For large problem instances this is overly  time consuming. Quasi-Newton methods have been developed that approximate the Hessian, however the convergence analysis is weaker than the full method~\citep{DenM77}. There are a lots of works in the field of distributed second order optimization such as~\cite{zhang2015disco},~\cite{smith2018cocoa},~\cite{wang2017giant} and~\cite{crane2019dingo}. In this work we take an alternative approach that retains the linear-quadratic convergence of Newton's method~\citep{PilW17} and adapt it to the distributed, which was first introduced in~\cite{bartan2020distributed}, but communication efficient setting. The key idea is to ``sketch'' the Hessian at each iteration.

\subsection{Contribution}
In this paper, we use the distributed iterative Hessian sketch algorithm (DIHS) which was introduced by~\cite{bartan2020distributed} to solve the large-scale system identification problems in a more scalable manner. Specifically, our contributions are:
\begin{itemize}
 	\item We give a new proof and a different convergence rate  from~\cite{bartan2020distributed, wang2017giant} for the DIHS algorithm.
 	\item  We provide a convergence guarantees for the DIHS algorithm with various sketching schemes on the matrix least square problems, not limited to Gaussian sketches mentioned in~\cite{bartan2020distributed}.
 \item We show that DIHS algorithm is consistent with the non-asymptotic sample complexity bound $O(\frac{1}{\sqrt{N}})$ for learning the Markov parameters.
\end{itemize}

\section{Background }

\emph{Sketching} has become a popular method for scientific computing workflows that deal with massive data sets; or more precisely, massive matrices~\citep{Drim16}. We develop an iterative, distributed sketching algorithm for solving system identification problems that are formulated as least-squares problems of the form~\eqref{eq:leastsq}. Consider the overdetermined least-squares problem with problem data $A\in \R^{n\times d}$, $b\in \R^n$, where $n \gg d$:
\begin{equation*}
	\underset{{x} \in \R^{d}}{\operatorname{minimize}}~\|Ax-b\|_{2}^{2},
\end{equation*}
and let $x^\star$ denote an optimal solution.\footnote{Note that problem~\eqref{eq:leastsq} is equivalent to this problem after vectorization.} Assuming $A$ is dense and has no discernible structure, factorization-based approaches solve the problem in $O(nd^2)$ arithmetic operations. The \emph{sketch-and-solve} approach  constructs a matrix $S$ of dimension $m\times n$ (where $m\ll n$) and solves
\begin{equation*}
	x^{\sharp} \in \underset{{x} \in \R^{d}}{\operatorname{arg min}}~\|S(Ax-b)\|_{2}^{2},
\end{equation*}
instead of the original problem. The driving idea is that if $m$ is small, then solving this problem is easier than solving the original problem. Amazingly, letting $S$ be a random matrix chosen from an appropriate distribution (to be defined later) will suffice. The matrix $SA$ is called the the sketch of $A$ and $S$ is the sketching matrix or embedding matrix. If $S$ is chosen as a subspace embedding of $\ra([A~b])$ then the action of $S$ preserves geometry and one can show that the residuals satisfy $\|Ax^\sharp -b \|_2\le (1+\epsilon)\|Ax^\star-b\|_2$, with high probability, where $\epsilon>0$ is the distortion of the embedding. In practice, even when the residuals are close, there is no guarantee that $\|x^\sharp - x^\star \|_2$ will be small. This problem will be alleviated by \emph{iterative-sketching} methods derived by~\cite{PilW16} described later. However, the sketch-and solve framework highlights several important points: How do we choose the random embedding (the matrix $S$)? How is the projection dimension $m$ chosen? Indeed, for the bound above to hold we require $m\sim d\log (d)/\epsilon ^2$ \citep{Sar06}. A geometric-type bound by~\cite{PilW15} showed that when entries of $S$ are sub-Gaussian, one requires $m \ge \frac{c}{\delta^2}\mathcal{W}(A\mathcal K)$ to obtain similar quality bounds in residual, where $\mathcal{W}(A\mathcal K)$ is the Gaussian width~\citep{vershynin2018high} of the cone $A\mathcal K$. Clearly there is tension between making $m$ small (improved computation) and obtaining an accurate solution.

 When $S$ is selected to be a dense matrix, the cost of forming $SA$ is $O(mnd)$ and no computational saving is achieved. However, there exist families of randomized matrices that admit a fast matrix-vector multiply which reduce the cost of forming the product to $O(nd\log(m))$  thus providing significant savings. In addition to the dense sub-Gaussian case, we will consider   of randomized embeddings defined by; randomized orthogonal systems (ROS)~\cite{AilC09}, Sparse Johnson-Lindenstrauss Transforms (SJLTs) \citep{KanN14}, and uniform sampling.

\section{Distributed Iterative Hessian Sketch}\label{sec:main}
Consider the optimization problem~\eqref{eq:leastsq}. Applying Newton's method  with a variable step-size $\alpha_t$, produces a sequence of iterates of the form:
\begin{equation}\label{eq:Newton}
X_{t+1} = X_t -\alpha_t(U^TU)^{-1}U^T(UX_t - Y), \quad t=1,2,\dots
\end{equation}
where  the Hessian is given by $U^TU$ and the gradient by $U^T(UX_t-Y)$ c.f.,~\citep{convex}. Given that the columns of $U$ are made up from the control input which we assume to be random Gaussian variables, $U$ will be a dense matrix. Under the assumption that $(N,m,T)$ are large forming the Hessian and the gradient will grind the Newton iterates to a halt.

Instead of computing the iterate~\eqref{eq:Newton} exactly, \cite{PilW16} introduced the iterative Hessian sketch (IHS) to approximate the Hessian. The explicit update rule is  given by:
\begin{equation}\label{eq:IHS}
X_{t+1} = X_t -\alpha_t(U^TS_t^T S_tU)^{-1}U^T(UX_t - Y),
\end{equation}
where $S_t \in \mathbb{R}^{s\times N} (s \ll N)$ denotes the embedding matrix at the  $t^{\text{th}}$ iteration. The matrix    $U^T S^T S U$ is called the the ``sketched Hessian''. We will now introduce some specific classes of sketching matrices.
\begin{definition}\label{def:subg}
A random variable $x$ such that $\mathbb E x=0$ is said to be sub-Gaussian with variance proxy $\sigma$, if its moment generating function satisfies
\[
\mathbb E\{ \exp(\lambda x)\}\le \exp\left({\frac{\sigma^2\lambda^2}{2}}\right).
\]
\end{definition}
An equivalent characterization   of a sub-Gaussian random variable obtained from Markov's inequality is that $\mathbb P(|x|\ge \lambda)\le 2\exp(-t^2/4)$, where $t=2\lambda$ and $\sigma=1$. We write $x\in \sg (\sigma^2)$ to denote that $x$ is sub-Gaussian. Note that this notation is not precise in the sense that $\sg(\sigma^2)$ denotes a family of distributions. The particular choice of distribution will be clear from context.
We consider the following families of random sketching matrices from which we draw  $S\in\R^{s\times N}$:
\begin{itemize}
\item \textbf{Sub-Gaussian:} Each element of $S$ is drawn from a specific sub-Gaussian distribution, i.e., $S_{ij}\iid \sg(\sigma^2) $ where the particular distribution is fixed for all entries. Examples of distributions that satisfy Definition~\ref{def:subg} include Gaussian, Bernoulli, and more generally, any bounded distribution.  Note that the sub-class of Gaussian sketch matrices are almost surely dense, thus they are often useful for proving results, and less useful  for computation. 
\item \textbf{Uniform}: Let $\{p_i\}_{i=1}^N$ denote the uniform distribution over $1,\hdots,N$. Then the uniform sketch samples  the rows $s$ times (with replacement). The $j^{\text{th}}$ row of $S$ is $s_j^T = e_j /\sqrt{p_j}$ with probability $p_j$, where $e_j$ is the $j^{\text{th}}$ standard basis vector. Other weights (probability distributions ) have been studied, however we do not pursue these here.
\item \textbf{Random Orthogonal System (ROS)-based Sketch:} This sketching matrix is based on a unitary trigonometric transform $F\in \mathbb F^{N \times N}$ (defined in  appendix~\ref{sec:Had}). As our least-squares problem is defined over the reals we restrict our attention to real transforms, and in particular the Walsh-Hadamard Transform. The matrix $S$ is then formed according to
\[
S = \sqrt{\frac{N}{s}}RFE,
\]
where $E=\mathrm{diag}(\nu_1,\dots , \nu_N)$ with $\nu_i$ drawn uniformly at random from $\{+1,-1\}$. The matrix $R$ is a $s\times N$ uniform sketching matrix defined above. The structure of an ROS matrix allows for a fast matrix-vector multiply.
\item \textbf{Sparse Johnson-Lindenstrauss Transform (SJLT)-based Sketches}. SJLT sketching matrices are another structured random matrix family that offer fast matrix-vector multiplication and are particularly suitable when the matrix to be sketched is sparse. Several constrictions exist, we follow that of \citep{KanN14}. Each column of $S$ has exactly $l$ non-zero entries at randomly chosen coordinates. The non-zero entries are chosen uniformly from $\{+1/\sqrt{l}, -1/\sqrt{l}\}$. SJLT matrices also belong to the class of sub-Gaussian sketching matrices.
\end{itemize}
Loosely speaking, ROS and SJLTs when applied to a vector attempt to evenly  mix all the coordinates and then randomly sample to obtain a lower dimensional vector with norm proportional to the original vector. In contrast, uniform sampling simply selects a subset of rows of $A$ chosen uniformly at random. A dense sub-Gaussian sketching matrix extends uniform sampling by linearly weighting the entries of each row.

Compared to $O((mT)^2 N)$ flops given by pseudo-inverse method, IHS with ROS or SJLT sketches takes $O((NmT \log(mT))\log(1/\epsilon))$ flops, which is linear in $NmT$, to obtain an $\epsilon$-accurate solution. Obviously, IHS has significantly lower complexity than direct method since we assume $N \gg mT$.

\begin{remark}
	IHS can also deal with the least square problem with $N\ll mT.$ In this case, we just need to sketch the column-space instead of the row-space. The constrained case is also easily handled.
\end{remark}

Just as with the Gaussian Newton Sketch, which produces unbiased estimates of the exact Newton step, many sketching matrices provide near-unbiased estimates of the Newton step~\citep{derezinski2021sparse}. This property is very important in distributed setting, where we can compute the iterate \eqref{eq:IHS} multiple times in parallel. Averageing schemes can then be employed to achieve better estimation performance~\citep{derezinski2019distributed,wang2017giant,wang2017sketched}. Using this idea, \cite{bartan2020distributed} introduced the Distributed-IHS (DIHS) algorithm, which is described by Algorithm~\ref{alg:DIHS}.

\begin{algorithm} 
	\caption{Distributed Iterative Hessian Sketch (DIHS)}
	\label{alg:DIHS}
	\begin{algorithmic}[1] 
		\State \textbf{Inputs}: Input matrix $U \in \mathbb{R}^{N\times mT}$, output matrix $Y \in \mathbb{R}^{N \times p}$, sketching size $s \ll N$.
		\State \textbf{Initialize}: Initial iterate $X_0 \in \mathbb{R}^{mT \times p}$
		\State \textbf{for} $t = 0,1,\cdots,M-1$ \textbf{do}
		\State \quad \textbf{Central node:} broadcasts $X_t$
		\State \quad \textbf{for} worker $i=1,2,\cdots,r$ \textbf{do in parallel}
		\State \quad \quad Generate a sketching matrix $S_{i}^t \in \mathbb{R}^{s \times N}$
		\State \quad \quad  Compute gradient $g_{t}=U^{T}\left(U X_{t}-Y\right)$.
		\State \quad \quad 
	$X_i^t = \underset{X}{\arg\min}  \Big\{  \frac{1}{2s} \lVert S_i^t U(X-X_t)\rVert_2^2 + \langle g_t, X\rangle \Big\}$
        \State \quad \quad Send $X_i^t$ to the central node
		\State \quad \textbf{end for}
		\State \quad \textbf{Central node:} Update $X_{t+1} = \frac{1}{r}\sum_{i=1}^r X_i^t$
		\State \textbf{end for}
	\end{algorithmic}
\end{algorithm}
At the $t^{\mathrm{th}}$ iteration of Algorithm~\ref{alg:DIHS}, each worker $i$ only has a small sketch of the full data set and computes the sketched version of Hessian matrix $U^T (S_i^t)^T S_i^t U$ and then computes the local update direction using the sketched Hessian. They then send the updated states to the central node. The central node averages all the states to update the new iteration $X_{t+1}$ and then broadcasts $X_{t+1}$ to all the workers. Using Algorithm~\ref{alg:DIHS},   the communication complexity decreases from $O((mT)^2)$ to $O(mT)$ at each iteration since we don't broadcast the Hessian to each worker, and only communicate the update direction.

\begin{remark}
	Note that the size of the sketching matrix can be different for each worker, i.e., $s\rightarrow s_i$ in line 6. This is particulary useful as it allows for the use of a heterogenous set of worker machines, each with their own resource profile.
\end{remark}

\subsection{Convergence Analysis}
We are now ready to state the main results of this work; the convergence analysis of the DIHS algorithm in terms of the approximation error and the estimation quality.
\begin{theorem}\label{thm1}
	Fix $\rho \in(0,1 / 2)$. If the number of rollouts satisfies $N \geq c T m \log ^{2}(2 T m) \log ^{2}(2 \bar{N} m)$ and if
	\begin{itemize}
		\item   $S_i^t$ is a sub-Gaussian sketching matrix with a sketching size $s \ge \frac{c_0}{\rho^2}mT,$ with probability at least $1-(2 \bar{N} m)^{-\log ^{2}(2 T m) \log (2 \bar{N} m)} -c_{1}  r t e^{-c_{2} s \rho^{2}}$,
		\end{itemize}
		or, 
		\begin{itemize}
		\item    $S_i^t$ is a randomized orthogonal system (ROS) sketching matrix with a sketching size $s \ge \frac{c_0 \log^4(mTp)}{\rho^2}mT,$ with probability at least $1-(2 \bar{N} m)^{-\log ^{2}(2 T m) \log (2 \bar{N} m)} -c_{1}  r t e^{-c_{2} \frac{s \rho^{2}}{\log ^{4}(mTp)}}$,
	\end{itemize}
    then, the output $X_t$ given by the DIHS algorithm at the $t^{\mathrm{th}}$ iteration satisfies
	$$
	\lVert X_t-X^{\mathrm{LS}}\rVert_{F}\leq 2\Big( \frac{\rho}{\sqrt{r}}\Big)^t \lVert X^{\mathrm{LS}}\rVert_{F},
	$$
where $X^{\mathrm{LS}}$ denotes the least-squares solution to problem~\eqref{eq:leastsq}, and $c, c_0, c_1$ and $c_2$ are absolute constants.
\end{theorem}
\begin{proof}
See appendix~\ref{sec:proofs} for this and subsequent proofs.
\end{proof}
\begin{remark}
The DIHS algorithm converges geometrically to the least-squares solution of~\eqref{eq:leastsq}. The linear convergence rate is $\frac{\rho}{\sqrt{r}}$, which decreases when number of workers $r$ increases. If we apply the DIHS algorithm with $O(\frac{\log(1/\epsilon)}{\log(\sqrt{r}/\rho)})$ iterations and choose the sketching matrices to satisfy the requirement of Theorem~\ref{thm1}, then the output which we denote by $\hat{X}$ satisfies:
	$$
\frac{\lVert \hat{X}-X^{\mathrm{LS}}\rVert_{F}}{ \lVert X^{\mathrm{LS}}\rVert_{F}}\leq \epsilon
$$	
with high probability.
\end{remark}

Next, we  quantify the approximation quality in terms of the distance between the DIHS solution and the ground-truth Markov parameters $G$ in the following theorem.

\begin{theorem}\label{thm2}
	Frame the hypotheses of Theorem~\ref{thm1}. For all $t\ge 1$, the output $X_t$ given by the DIHS algorithm satisfies:
\begin{equation}\label{eq:error}
	\lVert X_t - G^T \rVert_F  \le 2\Big( \frac{\rho}{\sqrt{r}}\Big)^t \lVert X^{\mathrm{LS}}\rVert_{F} + \frac{\left(\sigma_{v}+\sigma_{e}\right) \sqrt{p}+\sigma_{w}\|F\|_{2}}{\sigma_{u}} \sqrt{\frac{T q \log ^{2}(T q) \log ^{2}(N q)}{N}}
\end{equation}
with high probability.
\end{theorem}

\begin{remark}
	Note that the first term of the RHS of \eqref{eq:error} linearly converges to $0$ when $t \rightarrow \infty.$ Therefore, the estimation error $\lVert X_{t} - G^T \rVert_F$ given by the DIHS algorithm still maintains the $O(\frac{1}{\sqrt{N}})$ sample complexity when the number of iterations $t$ becomes large.	
\end{remark}

\section{Numerical Simulations}
\begin{figure}
	\centering     
\includegraphics[width=0.3\linewidth]{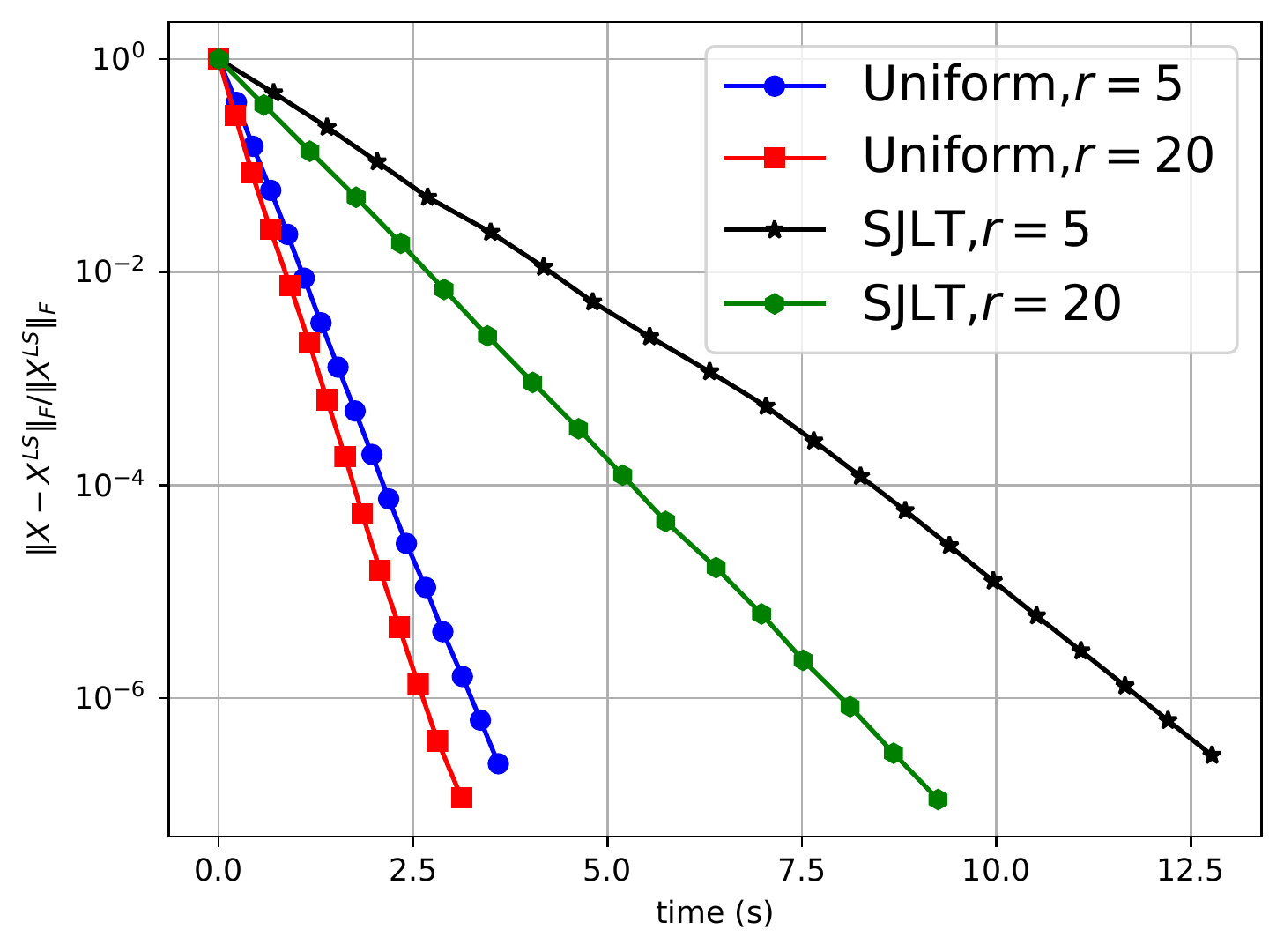}\quad
\includegraphics[width=0.3\linewidth]{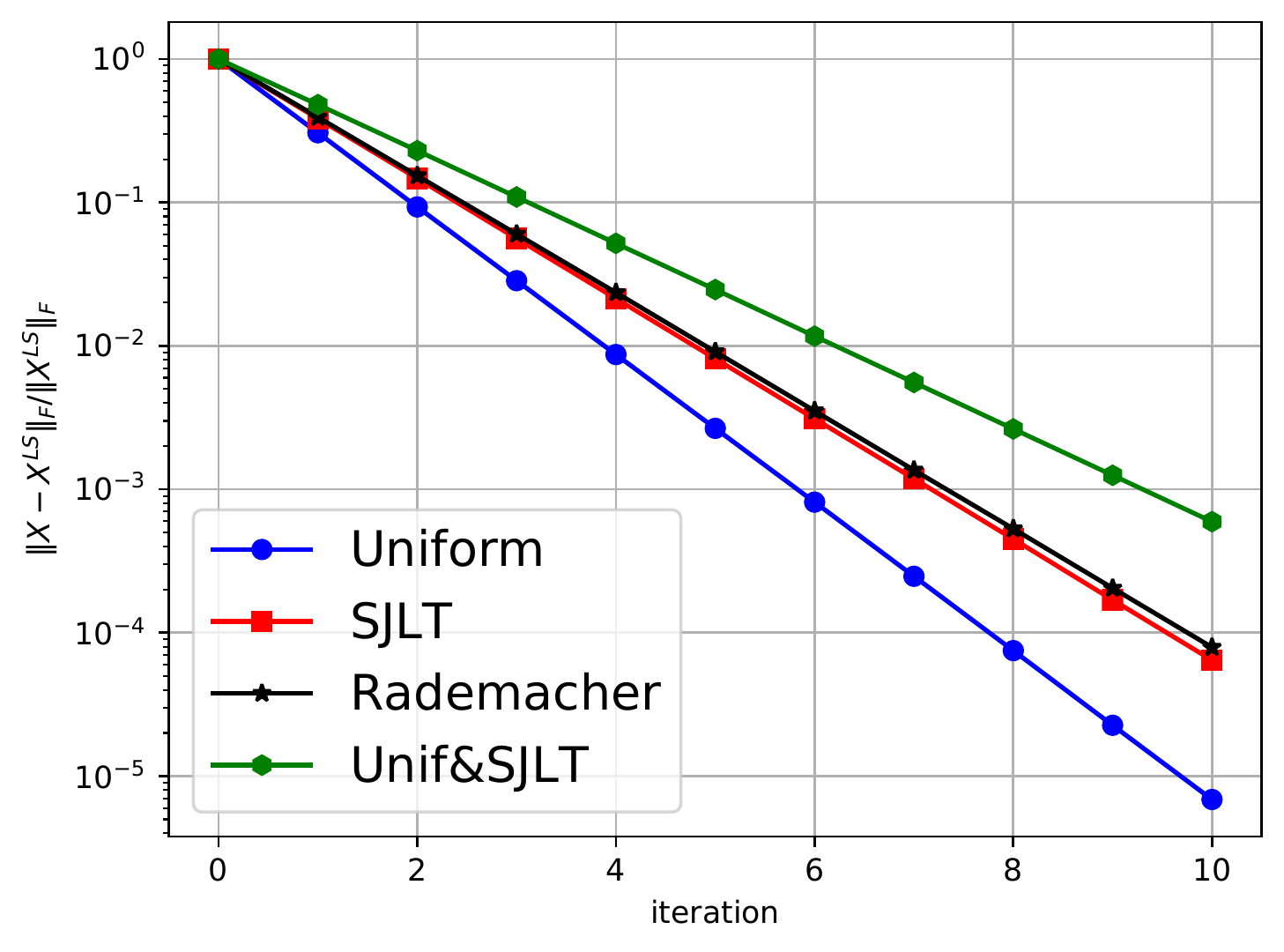}\quad
\includegraphics[width=0.3\linewidth]{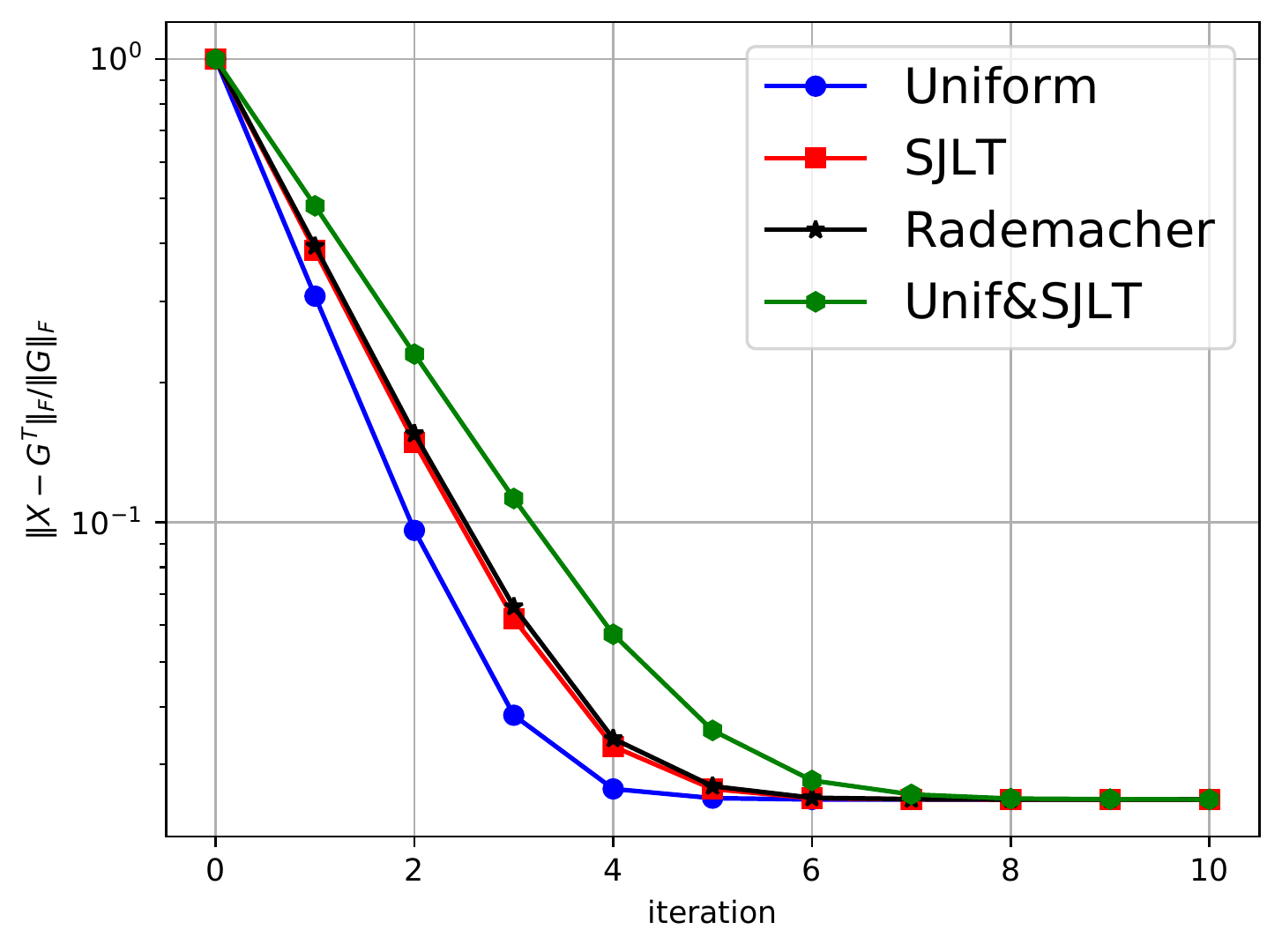}
	\caption{$(n,m,p,N,T)=(80,60,70,29971,30)$, $s=7200, s_1=14400$. The dimension of matrix $U$ is $(29971,1800)$ and matrix $Y$ is $(29971,70)$}
	\label{fig:a}
\end{figure}

\begin{figure}
	\centering     
	\includegraphics[width=0.3\linewidth]{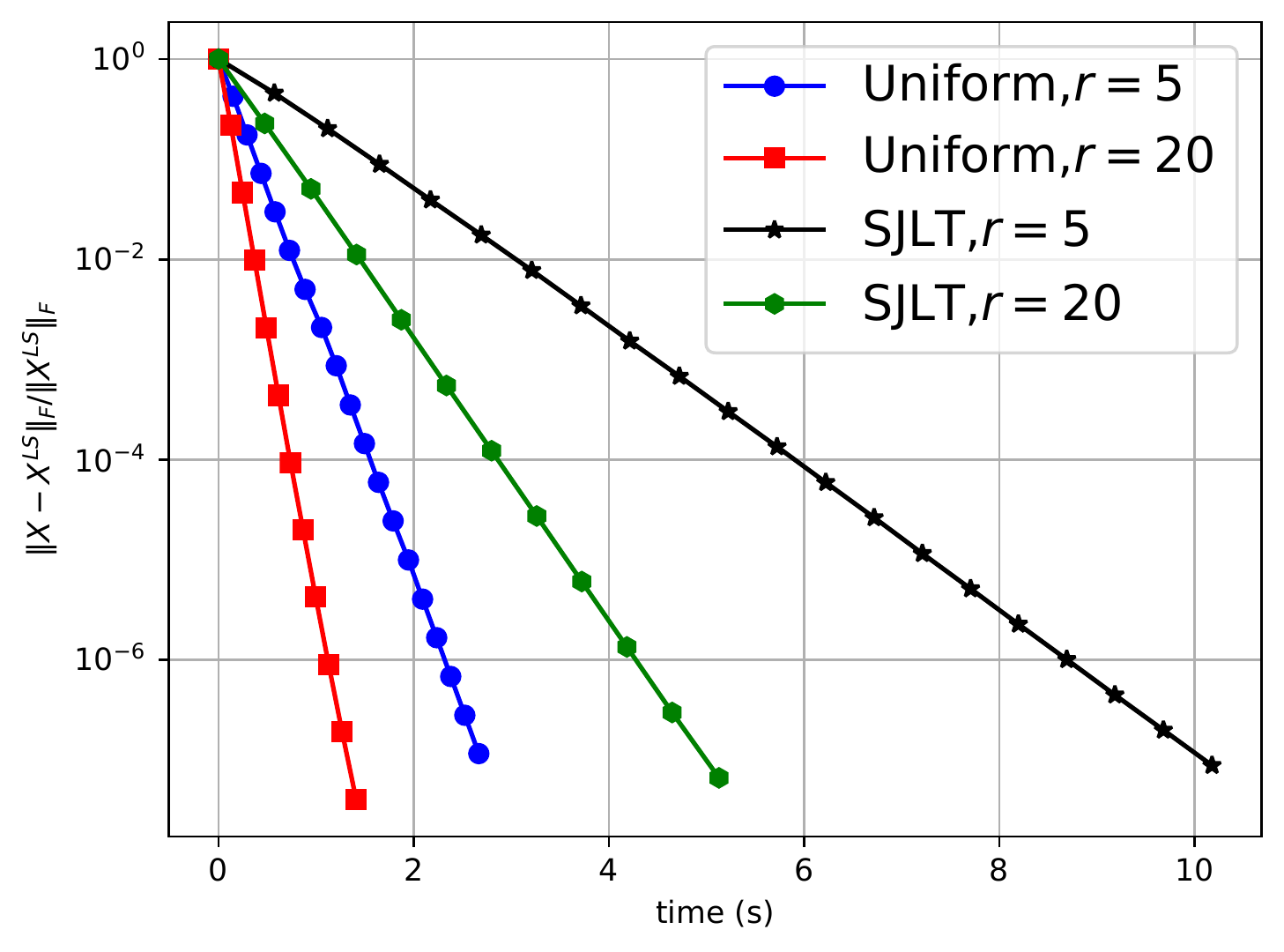}\quad
	\includegraphics[width=0.3\linewidth]{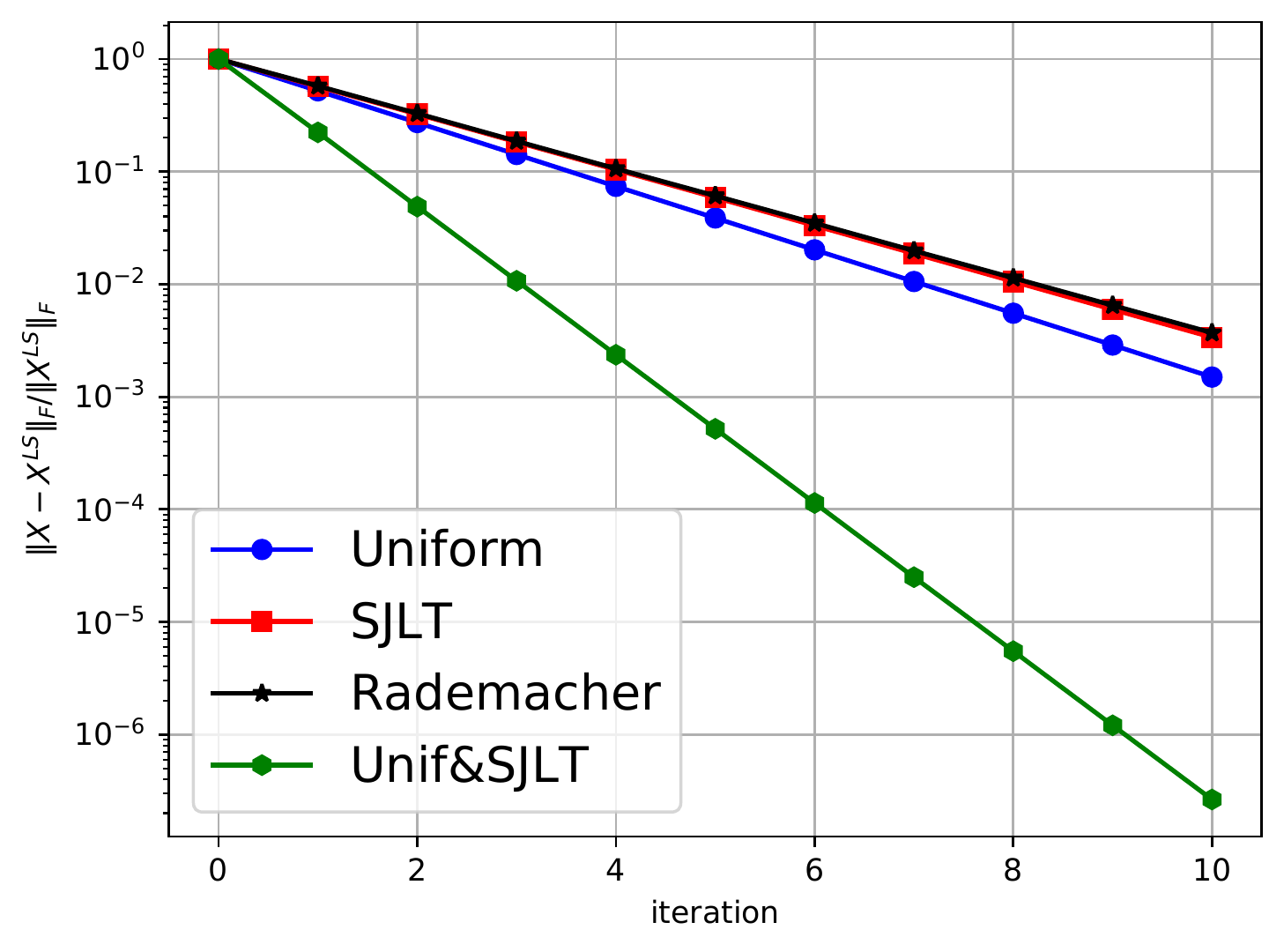}\quad
	\includegraphics[width=0.3\linewidth]{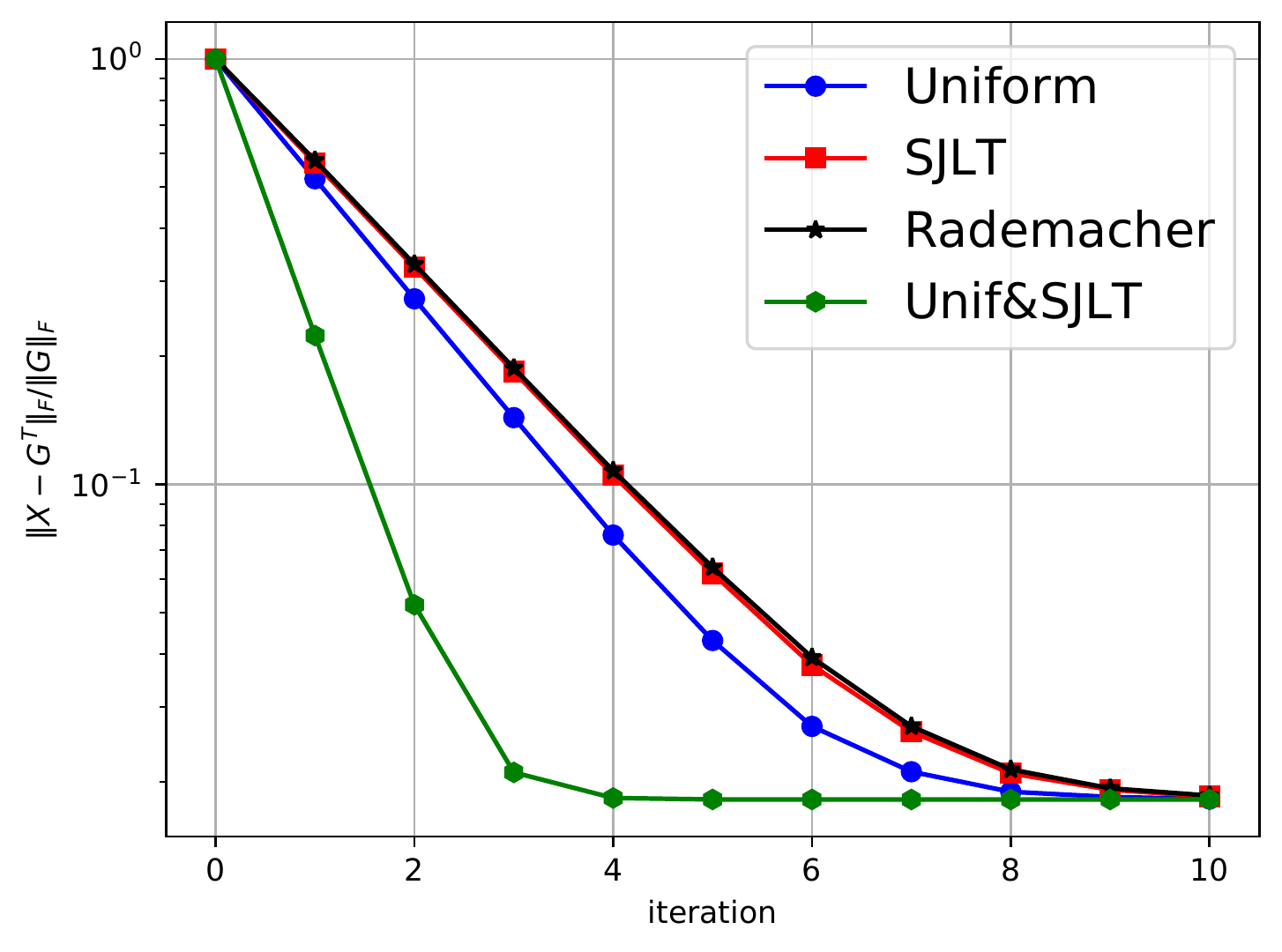}
	\caption{$(n,m,p,N,T)=(100,80,70,49981,20)$, $s=4800, s_1=7200$. The dimension of matrix $U$ is $(49981,1600)$ and matrix $Y$ is $(49981,70)$}
	\label{fig:b}
\end{figure}

\begin{figure}
	\centering     
	\includegraphics[width=0.3\linewidth]{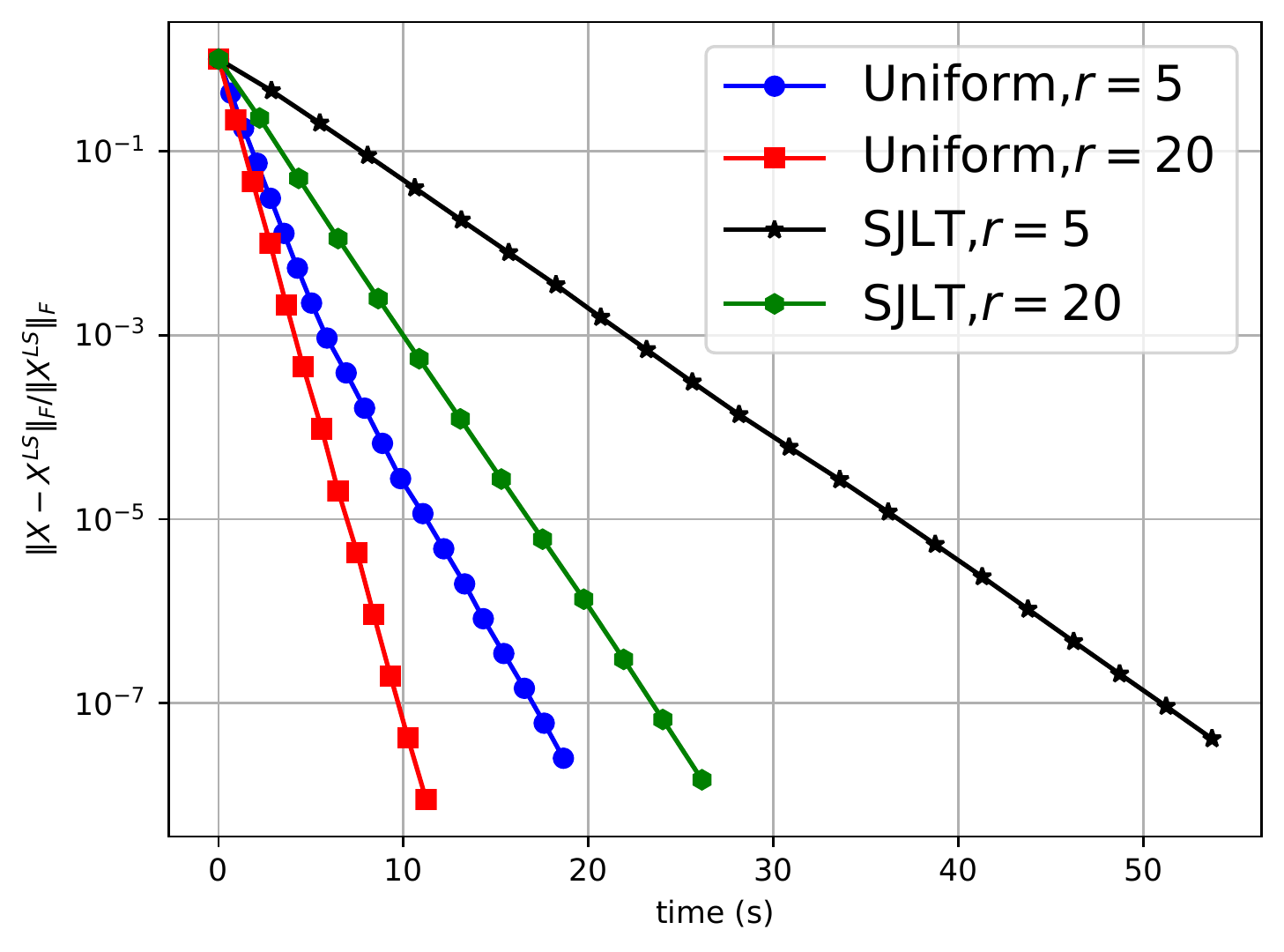}\quad
	\includegraphics[width=0.3\linewidth]{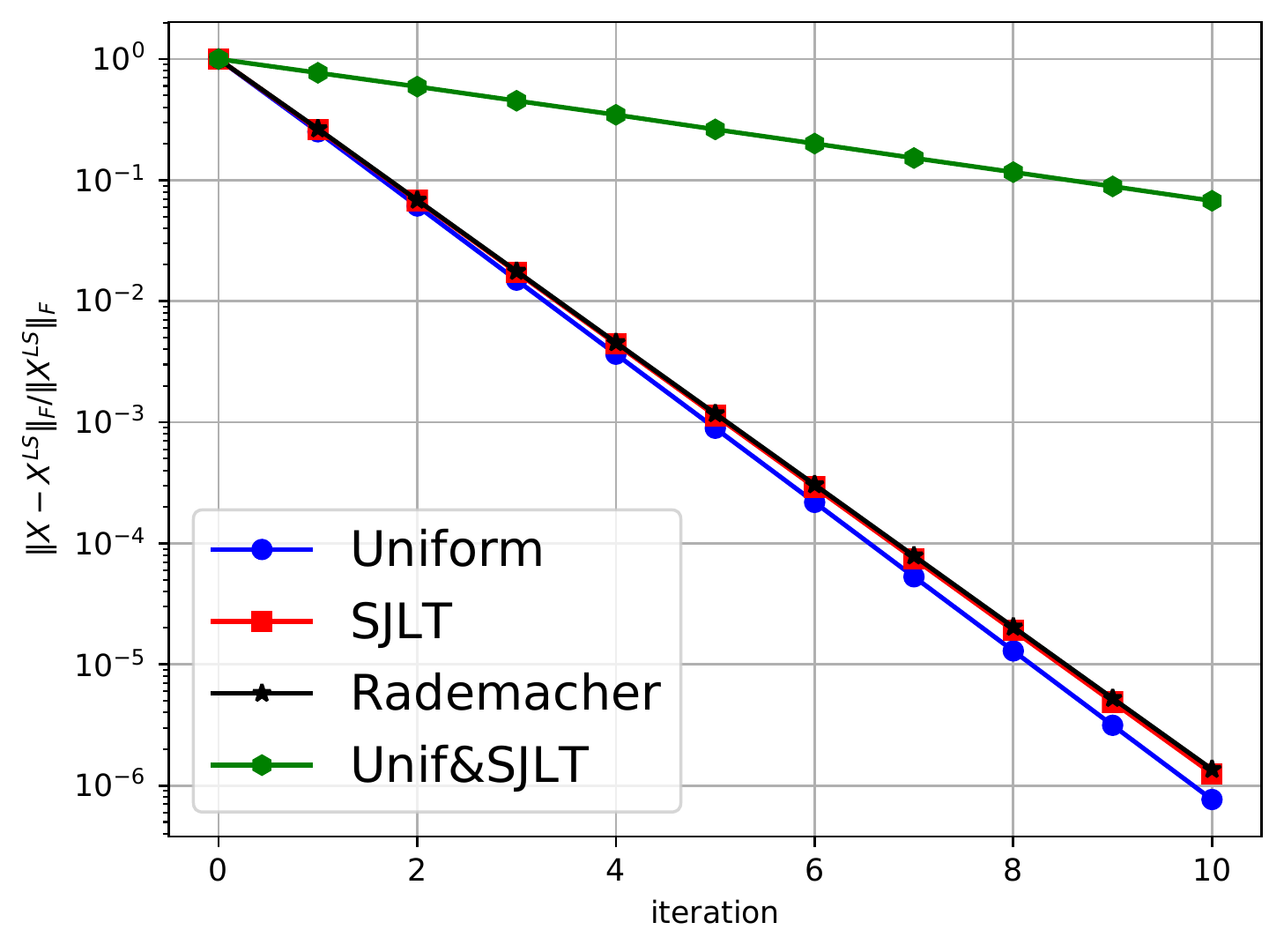}\quad
	\includegraphics[width=0.3\linewidth]{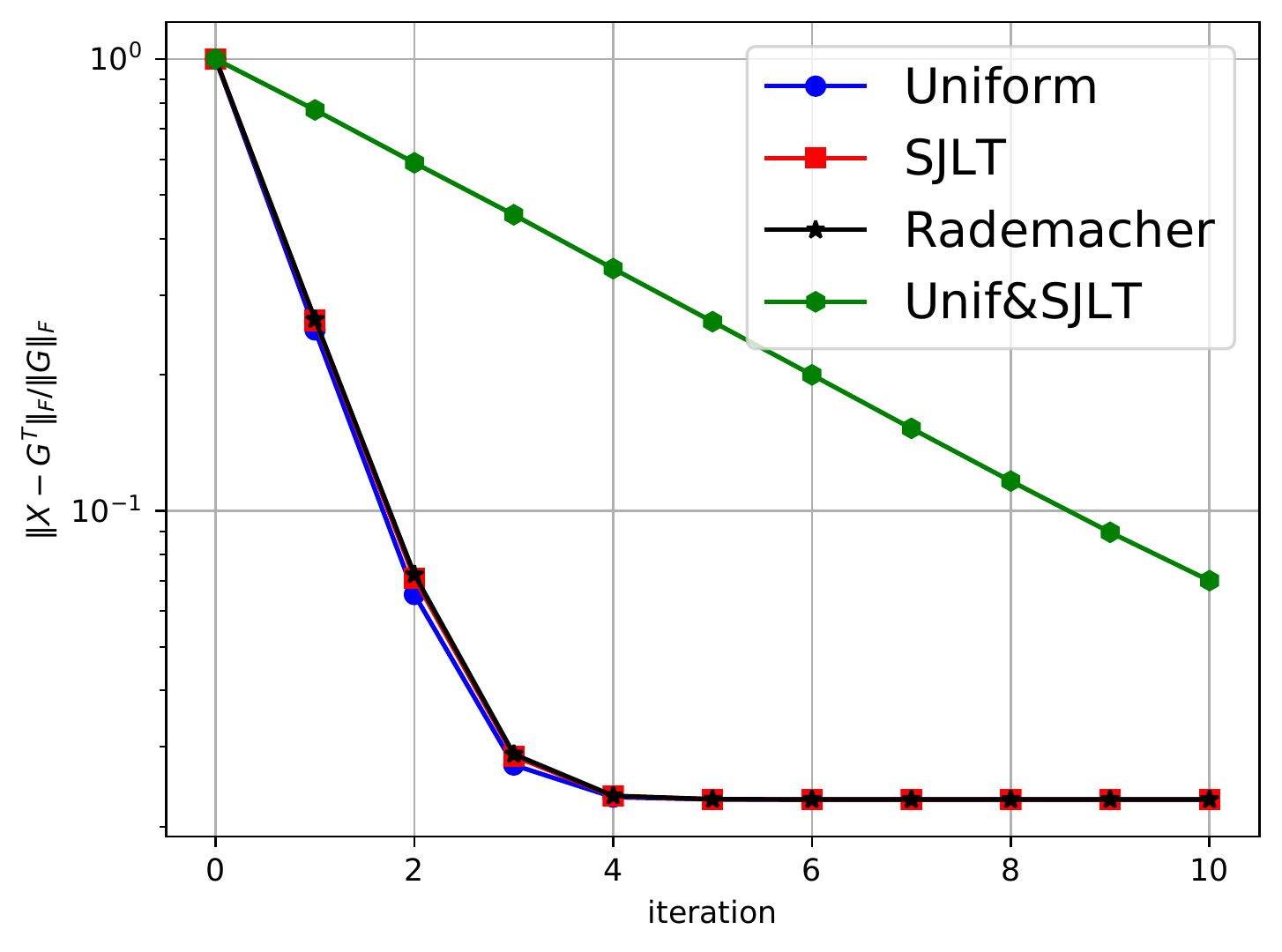}
	\caption{$(n,m,p,N,T)=(200,150,100,59981,20)$, $s=6000, s_1=7200$. The dimension of matrix $U$ is $(59981,3000)$ and matrix $Y$ is $(59981,100)$}
	\label{fig:c}
\end{figure}

We now demonstrate the performance of the DIHS algorithm on three randomly generated large-scale dynamic systems described by \eqref{eq:fully observed}. 
For each system, we choose the parameters $(n,m,p,N,T)$ as shown in the caption of Figures~\ref{fig:a}--\ref{fig:c}. To ensure a fair comparison, we fix a constant sketch dimension for all workers. Recall that theoretically this is unnecessary. We introduce two new sketching matrices. Rademacher sketches are defined such that each entry of $S$ is $\frac{1}{\sqrt{s}}$ with probability $1 / 2$ and $-\frac{1}{\sqrt{s}}$ otherwise. A two-stage uniform+SJLT sketch is produced from $S(S_1U)$ where $S$ is an SJLT sketching matrix with $s$ rows and $S_1$ is a uniform sketchimg matrix with $s_1$ rows.

We generate the system matrices $(A,B,C,D)$ through a uniform distribution over a range of integers as follows; entries of the matrix $A$ with random integers from $1$ to $5$, and matrices $B, C, D$ with random integers from $-2$ to $2$. Then, we re-scale the matrix $A$ to make it Schur-stable, i.e., $\lvert \lambda_{\mathrm{max}}(A)\rvert <1$. The standard deviations of the process and measurement noises are chosen to be $\sigma_w = 0.1$ and $\sigma_v= 0.1$. We fix the input variance at $\sigma_u=1$.

The left plot in each figure shows the normalized difference between the estimated solution at each iteration and the optimal least square solution $\|X^{\mathrm{LS}}\|$, (i.e. $\frac{\lVert X_t - X^{\mathrm{LS}}\rVert_F}{\lVert X^{\mathrm{LS}}\rVert_F}$) versus time (seconds) for the DIHS algorithm. In each system, we tested the performance of DIHS algorithm using the uniform and SJLT sketch with $r=5$ and $r=20$ worker machines. The  stopping criterion is that the distance between two consequent outputs (i.e. $\lVert X_{t+1} -X_{t} \rVert_F$) is less than $10^{-3}$. As predicted by the theoretical analysis, no matter what sketching matrix we use in the DIHS algorithm, the convergence rate decays as the number of workers $r$ increases. Compared to SJLT sketches, it seems that uniform sketching matrices could speed up the convergence throughout these three system identification examples. We note that in terms of computation speed, the uniform sketches should be fast as they require fewer arithmetic operations to apply and less time to construct than every other sketch type.  

The middle plot of  each figure illustrates the relative error between the estimated solution at the  $t^{\mathrm{th}}$ iteration and the optimal least square solution (i.e. $\frac{\lVert X_t - X^{\mathrm{LS}}\rVert_F}{\lVert X^{\mathrm{LS}}\rVert_F}$) against iteration for the DIHS algorithm with different sketching matrices for a fixed number of workers: $r=15.$ Finally, the  third column shows the relative error between the estimated solution and true Markov parameters (i.e. $\frac{\lVert X_t - G^{T}\rVert_F}{\lVert G\rVert_F}$) versus iteration with $15$ workers. 
From the figures, we can easily observe that DIHS algorithm with all these four sketching matrices converges geometrically to the least square solution, which is consistent with the analysis derived from Theorem~\ref{thm1}. For all sketching matrices we tested, DIHS algorithm can sucessfully learn the true Markov parameter $G$, which is stated in Theorem~\ref{thm2}. The performance of different sketching matrices depends on the choice of sketching size $s$ and $s_1$. 

\section{Conclusion}
We have demonstrated that a randomized version of Newton's algorithm can solve large-scale system identification problems and is consistent with recent state-of-the art sample complexity results. Geometric convergence was proven for all the standard sketching matrices and the dimension-dependence of the sketching matrix was also derived.  Future work will involve benchmarking this second-order method against distributed first-order methods such as accelerated stochastic gradient descent. We are currently integrating this work with our previous results which uses a randomized SVD to produce a system realization~\citep{WanA21},  with the goal of producing  end-to-end bounds.

\acks{HW is generously funded by a Wei family fellowship and the Columbia Data Science Institute. We also acknowledge funding from the DoE under grant DE-SC0022234.}

\bibliography{references,references_JA}
\appendix
 
	\section{Walsh-Hadamard Transform}\label{sec:Had}
The Walsh-Hadamard transform can be thought of as a generalized Fourier transform built from Hadamard matrices. A Hadamard matrix~\citep{HedW78} is an $n\times n$ matrix with mutually orthogonal columns whose elements take values $\pm 1$. As a result, for any Hadamard matrix $H\in \R^{n\times n}$, we have $H^TH=nI$. A recursive construction is a follows: define $H_0 = 1$, then 
\[
H_m=\frac{1}{\sqrt 2}\left[ 
\begin{array}{cc}
H_{m-1} & H_{m-1}\\
H_{m-1} & -H_{m-1}\end{array}
\right].
\]
To apply the Walsh-Hadamard transform to matrices that are not of compatible row dimension we simply pad with zeros.
 \section{Proofs and auxiliary results}\label{sec:proofs}
In this section, we prove the main results: Theorem~\ref{thm1} and~\ref{thm2}. We leverage the proof technique from~\cite{PilW16}. However the results need to be adapted to  target the matrix least-squares problem and to account for the distributed problem setting.

We use matrix $X \in \mathbb{R}^{mT \times p}$ to denote the possible descent matrix and vector $x \in \mathbb{R}^p$ to denote any column vector in the possible descent matrix $X$.  Here we assume a more general formulation where we append the convex constraint $X\in \mathcal C$ to Eq~\eqref{eq:leastsq}.
The set of possible descent directions $\left\{ x-x^{\mathrm{LS}} \mid x \in \mathcal{C} \right\}$ plays an important role in controlling the error w.r.t. the least square solution $X^{\mathrm{LS}}$. The affinely transformed tangent cone is defined as
\begin{equation}\label{eq:cone}
	\mathcal{K}^{\mathrm{LS}}=\left\{v \in \mathbb{R}^{N} \mid v=t U \left(x-x^{\mathrm{LS}}\right) \quad\right. \text{for some} \ t \geq 0 \ \text{and}  \left.x \in \mathcal{C}\right\}
\end{equation}

The error bound of DIHS algorithm largely depends on the following two quantities:
\begin{equation}
	\begin{aligned}
		Z_{1}(S):=& \inf _{v \in \mathcal{K}^{\mathrm{LS}} \cap \mathcal{S}^{N-1}} \frac{1}{s}\lVert S v\rVert_{2}^{2} \quad \text{and} \\
		Z_{2}(S):=& \sup _{v \in \mathcal{K}^{\mathrm{LS}} \cap \mathcal{S}^{N-1}}\left|\left\langle u,\left(\frac{S^{T} S}{s}-I_{n}\right) v\right\rangle\right|
	\end{aligned}
\end{equation}
where $u$ is any fixed unit-norm vector and $\mathcal{S}^{N-1}$ denotes the Euclidean sphere, i.e., $\mathcal{S}^{N-1}=\left\{z \in \mathbb{R}^{N} \mid\|z\|_{2}=1\right\}$.

Following~\cite{PilW16}, a ``good event" is defined as:
\begin{equation}
	\mathcal{E}(\rho):=\left\{Z_{1}(S) \geq 1-\rho, \text { and } Z_{2}(S) \leq \frac{\rho}{2}\right\} \quad 
\end{equation}
where $\rho \in\left(0, \frac{1}{2}\right)$ is a given tolerance parameter.

In the following lemma, we will describe how to choose the sketching size $s$ in order to ensure the good events $\mathcal{E}(\rho)$ holds with high probability. Before the statement of lemma, we need to introduce the notion of Gaussian width:
$$
\mathcal{W}\left(\mathcal{K}^{\mathrm{LS}}\right):=\mathbb{E}_{g}\left[\sup _{v \in \mathcal{K}^{\mathrm{LS}}\cap \mathcal{S}^{N-1}}|\langle g, v\rangle|\right]
$$
where $g \sim N\left(0, I_{n}\right)$ is a standard Gaussian vector. 

\begin{remark}
	 It is easy to show that the Gaussian width of unit sphere $\mathcal{S}^{N-1}$ is at most $\sqrt{N}$. And from~\cite{PilW15}, we have that $\mathcal{W}\left(\mathcal{K}^{\mathrm{LS}}\right) \le \sqrt{mT}.$
\end{remark}

\begin{lemma}\label{sample}(Sufficient conditions on sketch dimension~\cite{PilW15})
	\begin{enumerate} 
\item For sub-Gaussian sketch matrices, given a sketch size $s>\frac{c_{0}}{\rho^{2}} \mathcal{W}^{2}\left(\mathcal{K}^{L S}\right)$, we have
		$$
		\mathbb{P}[\mathcal{E}(\rho)] \geq 1-c_{1} e^{-c_{2} s \delta^{2}}
		$$
		\item For randomized orthogonal system (ROS) sketches (sampled with replacement), given a sketch size $s>\frac{c_{0} \log ^{4}(mTp)}{\rho^{2}} \mathcal{W}^{2}\left(\mathcal{K}^{L S}\right)$, we have
		$$
		\mathbb{P}[\mathcal{E}(\rho)] \geq 1-c_{1} e^{-c_{2} \frac{s \rho^{2}}{\log ^{4}(mTp)}}
		$$
	\end{enumerate}
\end{lemma}

In addition, we define the sequence of ``good events"
\begin{equation}
	\mathcal{E}_i^t(\rho):=\left\{Z_{1}(S_{i}^t) \geq 1-\rho, \text { and } Z_{2}(S_{i}^t) \leq \frac{\rho}{2}\right\} \quad \text{for}\ \ i =1, \cdots, r \quad \text{and} \ \ t =0, 1,\cdots, M-1.
\end{equation}
Then we have the following error bound:
\begin{proposition}\label{pro4}
	For a fixed $\rho \in\left(0, \frac{1}{2}\right)$, the final solution $\hat{X} = X_M$ given by DIHS algorithm satisfies the bound
	\begin{equation}\label{eq:lea}
		\lVert \hat{X} - X^{\mathrm{LS}}\rVert_{U} \leq (\frac{\rho}{\sqrt{r}})^M \lVert X^{\mathrm{LS}}\rVert_{U},
\end{equation}
	conditioned on the event $\cap_{t=0}^{M-1} \cap_{i=1}^{r} \mathcal{E}^{t}_i(\rho)$. In \eqref{eq:lea}, $\lVert \cdot \rVert_{U}$ denotes $\lVert U \cdot \rVert_F.$ 
\end{proposition}
\begin{proof}
	If we can show that, for each iteration $t=0,1,2, \ldots M-1$, we have
	\begin{equation}\label{eq:success}
		\left\|X_{t+1}-X^{\mathrm{LS}}\right\|_{U} \leq \frac{\rho}{\sqrt{r}}\left\|X_{t}-X^{\mathrm{LS}}\right\|_{U},
	\end{equation}
	The claimed bounds \eqref{eq:lea} then hold by applying the bound \eqref{eq:success} iteratively from step 0 to $M-1$.
	
 Define $\Delta=X_{t+1}-X^{\mathrm{LS}}$. With some simple algebra calculation, we can rewrite $X_{t+1}$ as:
\begin{equation}\label{eq:opt}
X_{t+1}=\frac{1}{r} \sum_{i=1}^r \arg \min _{X}\left\{\frac{1}{2 s}\|S_{i}^t U X\|_{2}^{2}-\langle U^T \widetilde{Y}_i, X \rangle\right\},
\end{equation}
where $\widetilde{Y}_i^t:=Y-\left[I-\frac{(S_{i}^t)^{T} S_{i}^t}{s}\right] U X_{t}$. Since $X_{t+1}$ is the optimal solution for~\eqref{eq:opt}, the first-order optimality conditions tell us
\begin{equation}\label{eq:step1}
	\frac{1}{r}\sum_{i=1}^r \langle U^{T} \frac{(S_{i}^t)^{T} S_{i}^t}{s}U X_{t+1}-U^{T} \widetilde{Y}_i^t, X^{\mathrm{LS}}-X_{t+1}\rangle \geq 0
\end{equation}
Also, $X^{\mathrm{LS}}$ is optimal for the original least square program, we have
\begin{equation}\label{eq:step2}
	\frac{1}{r}\sum_{i=1}^r \langle U^{T}\left(UX^{\mathrm{LS}}-\widetilde{Y}_i^t-\left[I-\frac{(S_{i}^t)^{T} S_{i}^t}{s}\right] UX_{t}\right), X_{t+1}-X^{\mathrm{LS}}\rangle \geq 0
\end{equation}
Adding these two inequalities together, we get:
\begin{equation}\label{eq:internal}
	\frac{1}{r}\sum_{i=1}^r \frac{1}{s}\|S_i^t U \Delta\|_{F}^{2} \leq \frac{1}{r}\sum_{i=1}^r \langle U \left(X^{\mathrm{LS}}-X_{t}\right), \left(I-\frac{(S_{i}^t)^{T} S_{i}^t}{s}\right) U\Delta \rangle.
\end{equation}
According to the definition of $Z_{2}\left(S_{i}^t\right)$, we have
\begin{equation}\label{eq:step3}
	\begin{aligned}
	\langle U \left(X^{\mathrm{LS}}-X_{t}\right),& \left(I-\frac{(S_{i}^t)^{T} S_{i}^t}{s}\right) U\Delta \rangle = \sum_{j=1}^{p}\langle U\left(X^{\mathrm{LS}}(:,j)-X_{t}(:,j)\right), \left(I-\frac{(S_{i}^t)^{T} S_{i}^t}{s}\right) U \Delta(:,j) \rangle\\
		&\leq \sum_{j=1}^{p} \lVert U\left(X^{\mathrm{LS}}(:,j)-X_{t}(:,j)\right) \rVert_2 \lVert U \Delta(:,j)\rVert_2 Z_{2}\left(S_i^t\right)\\
		&\leq \left( \sum_{j=1}^{p} \lVert U\left(X^{\mathrm{LS}}(:,j)-X_{t}(:,j)\right) \rVert_2^2 \right)^{1/2}\left(\sum_{j=1}^p
		\lVert U \Delta(:,j)\rVert_2^2\right)^{1/2}Z_{2}\left(S_i^t\right)\\
		& \leq\frac{\rho}{2} \lVert U\left(X^{\mathrm{LS}}-X_{t}\right)\rVert_F \lVert U \Delta\rVert_F,
	\end{aligned}
\end{equation}
as long as the good event $\cap_{t=0}^{M-1} \cap_{i=1}^{r} \mathcal{E}^{t}_i(\rho)$ happens. In~\eqref{eq:step3}, $X^{\mathrm{LS}}(:,j)$ denotes $j$-th column of matrix $X^{\mathrm{LS}}.$ Define $$S_t=\frac{1}{\sqrt{r}}\left[(S_{1}^t)^{T}, \cdots,(S_{r}^t)^{T}\right]^T \in \mathbb{R}^{sr \times N}.$$ According to Theorem 14 in~\cite{wang2017sketched}, if 
\begin{equation}\label{eq:MMP}
	\frac{1}{r}\sum_{i=1}^r\langle U \left(X^{\mathrm{LS}}-X_{t}\right), \left(I-\frac{(S_{i}^t)^{T} S_{i}^t}{s}\right) U\Delta \rangle \leq\frac{\rho}{2} \lVert U\left(X^{\mathrm{LS}}-X_{t}\right)\rVert_F \lVert U \Delta\rVert_F
\end{equation}
holds with probability $p=\mathbb{P}(\cap_{t=0}^{M-1} \cap_{i=1}^{r} \mathcal{E}^{t}_i(\rho))$, then 
\begin{equation}\label{eq:TMMP}
\langle U \left(X^{\mathrm{LS}}-X_{t}\right), \left(I-\frac{(S_t)^{T} S_t}{s}\right) U\Delta \rangle \leq\frac{\rho}{2\sqrt{r}} \lVert U\left(X^{\mathrm{LS}}-X_{t}\right)\rVert_F \lVert U \Delta\rVert_F
\end{equation}
holds with the same probability. 

By the definition of $Z_{1}\left(S_{i}^t\right)$, the lefthand side of~\eqref{eq:internal} satisfies:
\begin{equation}\label{eq:RHS}
	\frac{1}{r}\sum_{i=1}^r \frac{1}{s}\|S_i^t U \Delta\|_{F}^{2} \ge Z_{1}\left(S_{i}^t\right) \lVert U \Delta\rVert_F^2 \ge (1-\rho)\lVert U \Delta\rVert_F^2
\end{equation}

Combining the inequality~\eqref{eq:RHS} and~\eqref{eq:TMMP} yields, as long as the good event $\cap_{t=0}^{M-1} \cap_{i=1}^{r} \mathcal{E}^{t}_i(\rho)$ occurs,
\begin{equation}
	\lVert	U \Delta\rVert_{F} \le \frac{\rho}{2(1-\rho)\sqrt{r}}\lVert U\left(X_{t}-X^{\mathrm{LS}}\right)\rVert_F \le \frac{\rho}{\sqrt{r}}\lVert U\left(X_{t}-X^{\mathrm{LS}}\right)\rVert_F,
\end{equation}
holds, then the proof is done.
\end{proof}

The following lemma tells us the lower bound and upper bound of singular values of matrix $U.$
\begin{lemma}\label{condi}\citep{oymak2019non}
	 Suppose the sample size $N \geq c T m \log ^{2}(2 T m) \log ^{2}(2 \bar{N} m)$ for sufficiently large constant $c>0 .$ Then, with probability at least $1-(2 \bar{N} m)^{-\log ^{2}\left(2 Tm\right) \log (2 \bar{N} m)}$, the input data matrix $U \in \mathbb{R}^{N \times T m}$ defined in Eq~\eqref{eq:leastsq} satisfies
	$$
	2 N \sigma_{u}^{2}I \geq U^*U \geq IN \sigma_{u}^{2} / 2.
	$$
\end{lemma}

\subsection{Proof of Theorem~\ref{thm1}}
Now we are ready to prove Theorem~\ref{thm1}. From Proposition~\ref{pro4}, we have
	\begin{equation*}
	\lVert \hat{X} - X^{\mathrm{LS}}\rVert_{U} \leq (\frac{\rho}{\sqrt{r}})^M \lVert X^{\mathrm{LS}}\rVert_{U},
\end{equation*}
	conditioned on the event $\cap_{t=0}^{M-1} \cap_{i=1}^{r} \mathcal{E}^{t}_i(\rho)$ for some $\rho \in (0,\frac{1}{2}).$ According to Lemma~\ref{condi}, $\kappa (U) \le 2$ holds with probability at least $1-(2 \bar{N} m)^{-\log ^{2}\left(2 Tm\right) \log (2 \bar{N} m)}$, as long as the number of rollouts $N \geq c T m \log ^{2}(2 T m) \log ^{2}(2 \bar{N} m)$, where $\kappa (U)$ is the condition number of matrix $U$. With these two facts, we can conclude that when $N \geq c T m \log ^{2}(2 T m) \log ^{2}(2 \bar{N} m)$ and the event $\cap_{t=0}^{M-1} \cap_{i=1}^{r} \mathcal{E}^{t}_i(\rho)$ happens, we have 
		\begin{equation*}
		\lVert \hat{X} - X^{\mathrm{LS}}\rVert_{F} \leq 2(\frac{\rho}{\sqrt{r}})^M \lVert X^{\mathrm{LS}}\rVert_{F},
	\end{equation*}
with probability at least $1-(2 \bar{N} m)^{-\log ^{2}\left(2 Tm\right) \log (2 \bar{N} m)}$.

According to lemma~\ref{sample}, if we apply the sub-Gaussian sketch matrices with sketching size $s>\frac{c_{0}}{\rho^{2}} mT$, we have	$$
\mathbb{P}[\mathcal{E}(\rho)] \geq 1-c_{1} e^{-c_{2} s \rho^{2}}.
$$
Applying the union bound, we  conclude that as long as $s>\frac{c_{0}}{\rho^{2}} mT$, then
$$
\mathbb{P}[\cap_{t=0}^{M-1} \cap_{i=1}^{r} \mathcal{E}^{t}_i(\rho)] \geq 1-c_{1}r M  e^{-c_{2} s \rho^{2}}.
$$
Using the union bound argument again, we have
		\begin{equation*}
	\lVert \hat{X} - X^{\mathrm{LS}}\rVert_{F} \leq 2(\frac{\rho}{\sqrt{r}})^M \lVert X^{\mathrm{LS}}\rVert_{F},
\end{equation*} holds with probability $1-(2 \bar{N} m)^{-\log ^{2}\left(2 Tm\right) \log (2 \bar{N} m)}- c_{1}r M  e^{-c_{2} s \rho^{2}}$ as long as the number of rollouts $N \geq c T m \log ^{2}(2 T m) \log ^{2}(2 \bar{N} m)$ and $s>\frac{c_{0}}{\rho^{2}} mT$. The proof is the same for ROS sketching matrices.

\subsection{Proof of Theorem~\ref{thm2}}
By using triangle inequality, we have
\begin{equation}\label{eq:pthm2}
	\begin{aligned}
		\|X_t-G^T\|_F &\leq \|X_t - X^{\mathrm{LS}}\|_F + \|X^{\mathrm{LS}} - G^T\|_F\\
		& = \|X_t - X^{\mathrm{LS}}\|_F +\|\hat{G}- G\|_F\\
		&\leq 2(\frac{\rho}{\sqrt{r}})^t \lVert X^{\mathrm{LS}}\rVert_{F} + \frac{\left(\sigma_{v}+\sigma_{e}\right) \sqrt{p}+\sigma_{w}\|F\|_{2}}{\sigma_{u}} \sqrt{\frac{T q \log ^{2}(T q) \log ^{2}(N q)}{N}}
		\end{aligned}
\end{equation}
holds with probability $1-(2 \bar{N} m)^{-\log ^{2}\left(2 Tm\right) \log (2 \bar{N} m)}- c_{1}r t  e^{-c_{2} s \rho^{2}}$ as long as the number of rollouts $N \geq c T m \log ^{2}(2 T m) \log ^{2}(2 \bar{N} m)$ and $s>\frac{c_{0}}{\rho^{2}} mT$. The last inequality of Eq~\eqref{eq:pthm2} is from Theorem~\ref{thm1} and \cite{oymak2019non}.

\end{document}